\documentclass{amsart}
\usepackage{amsmath,amsfonts}
\usepackage{latexsym}
\title[$\Psi$dos with symbols valued on a noncommutative algebra]
{Cordes characterization for pseudodifferential operators with symbols valued on a noncommutative C$^*$-algebra}
\author[S. T. Melo \and M. I. Merklen]{Severino T. Melo 
\and Marcela I. Merklen}
\date{}
\newtheorem{thm}{Theorem}
\newtheorem{pro}{Proposition}

\newtheorem{df}{Definition}
\begin{document}
\newcommand{\A}{{mathcal A}}
\newcommand{\R}{{\mathbb R}}
\newcommand{\C}{{\mathbb C}}
\newcommand{\N}{{\mathbb N}}
\newcommand{\um}{{\mathbf 1}}
\newcommand{\rn}{{\mathbb R}^{n}}
\newcommand{\rtn}{{\mathbb R}^{2n}}
\newcommand{\op}{operator}
\newcommand{\ops}{operators}
\newcommand{\psd}{pseudo\-dif\-fer\-en\-tial}
\newcommand{\Psd}{Pseudo\-dif\-fer\-en\-tial}
\newcommand{\bC}{{\mathcal B}^{\mathcal A}(\rtn)}
\newcommand{\bCk}{{\mathcal B}^{M_k(\mathcal A)}(\rtn)}
\newcommand{\bcn}{{\mathcal B}^{\mathcal A}(\rn)}
\newcommand{\sC}{{\mathcal S}^{\mathcal A}(\rn)}
\newcommand{\cqd}{\hfill$\Box$}
\newcommand{\h}{{\mathcal H}}
\newcommand{\pf}{{\em Proof}: }

\begin{abstract}
Given a separable unital C$^*$-algebra ${\mathcal A}$ with norm $||\cdot||$, let $E$ 
denote the Banach-space completion of the ${\mathcal A}$-valued Schwartz space on $\rn$ 
with norm $||f||_2=||\langle f,f\rangle||^{1/2}$, 
$\langle f,g\rangle=\int f(x)^*g(x)\,dx$. The assignment of the 
pseudodifferential operator $B=b(x,D)$ with ${\mathcal A}$-valued symbol 
$b(x,\xi)$ to each smooth function with bounded derivatives $b\in\bC$ 
defines an injective mapping $O$, from $\bC$ to the set $\h$ of all 
operators with smooth orbit under the canonical action of the Heisenberg group 
on the algebra of all adjointable operators on the Hilbert module $E$. It is known that
$O$ is surjective if ${\mathcal A}$ is commutative. In this paper, we show that, if $O$ 
is surjective for ${\mathcal A}$, then it is also surjective for $M_k({\mathcal A})$.
\end{abstract}
\maketitle

\begin{center}
\end{center}

\section{Introduction}

Let ${\mathcal A}$ be a separable C$^*$-algebra with norm $||\cdot||$ and unit $\um$, and
let $\sC$ denote the set of all ${\mathcal A}$-valued smooth (Schwartz) functions on $\rn$ which, 
together with all their derivatives, are bounded by arbitrary negative 
powers of $|x|$, $x\in\rn$. We equip it with the ${\mathcal A}$-valued inner-product 
\[
\langle f,g\rangle\,=\,\int f(x)^*g(x)\,dx,
\]
which induces the norm $||f||_2=||\langle f,f\rangle||^{1/2}$, and denote by 
$E$ its Banach-space completion with this norm. The inner product 
$\langle\cdot,\cdot\rangle$ 
turns $E$ into a Hilbert module \cite{L}. The set of all (bounded) adjointable \ops\ 
on $E$ is denoted ${\mathcal B}^*(E)$. 

Let $\bC$ denote the set of all 
smooth bounded functions from $\rtn$ to ${\mathcal A}$ whose derivatives of arbitrary 
order are also bounded. For each $b$ in $\bC$, a linear mapping from $\sC$ to 
itself is defined by the formula
\begin{equation}
\label{psddef}
(Bu)(x)\,=\,\frac{1}{(2\pi)^{n/2}}\int e^{ix\cdot\xi}b(x,\xi)\hat u(\xi)\,d\xi,
\end{equation}
where $\hat u$ denotes the Fourier transform, 
\[
\hat u(\xi)=(2\pi)^{-n/2}\int e^{-iy\cdot\xi}u(y)\,dy.
\] 
The operator $B:=b(x,D)$ extends to an element of 
${\mathcal B}^*(E)$ whose norm satisfies the following estimate. There 
exists $K>0$ depending only on $n$ such that 
\begin{equation}
\label{cvm}
||B||\,\leq\,K \sup\{\,
||\partial_x^\alpha\partial_\xi^\beta b(x,\xi)||;
(x,\xi)\in\rtn\ \mbox{and}\ \alpha,\beta
\leq(1,\cdots,1)\,\}.
\end{equation}
This generalization of the Calder\'on-Vaillancourt Theorem 
\cite{CV} was proven by Mer\-klen \cite{M,Mt}, see also \cite{H,R,S}.

The estimate \eqref{cvm} implies that the mapping 
\begin{equation}
\label{azz}
\rtn\ni(z,\zeta)\longmapsto  B_{z,\zeta}= T_{-z}M_{-\zeta}BM_{\zeta}
T_z\in{\mathcal B}^*(E)
\end{equation}
is smooth (i.e., $C^\infty$ with respect to the norm topology), where $T_z$
and $M_\zeta$ are defined by $T_zu(x)=u(x-z)$ and $M_{\zeta}u(x)=e^{i\zeta \cdot x}u(x)$, $u\in\sC$. 
That follows just as in the scalar case \cite[Chapter 8]{book}. 

\begin{df}
We call 
{\em Heisenberg smooth} an operator $B\in{\mathcal B}^*(E)$ for which the 
mapping \eqref{azz} is smooth, and denote by $\h$ the set of all such \ops. 
\end{df}
The elements of $\h$ are the smooth vectors for the canonical action of the Heisenberg 
group on ${\mathcal B}^*(E)$.

We therefore have a mapping
\begin{equation}
\label{o}
\begin{array}{rcl} 
O_{{\mathcal A}}:\bC&\longrightarrow&\h\\
b\ &\longmapsto&b(x,D).
\end{array}
\end{equation}
It is a standard result that, in the scalar case (${\mathcal A}=\C$), 
$O_{{\mathcal A}}$ is injective. For general ${\mathcal A}$, injectiveness follows from the scalar case by a 
duality argument. Cordes \cite{C}\ proved that $O_{{\mathcal A}}$ is surjective in the scalar case.
We have  shown \cite{MM2} that this also happens if ${\mathcal A}$ is unital and commutative. 

In this paper, we show that, if $O_{{\mathcal A}}$ is surjective, then $O_{M_k({\mathcal A})}$ is also surjective. 
We show that by first noticing that the Hilbert module $E_k$ for the matrix case is a 
Banach-space direct sum of $k^2$ copies of $E$. Then it follows that a bounded operator on $E_k$ (regarded only 
as a Banach space) is smooth under the action of the Heisenberg group if and only if it is 
a matrix whose entries are \ops\ on $E$ which are also smooth under the Heisenberg group. 
When we impose that such a matrix be an adjointable
$M_k({\mathcal A})$-module homomorphism, then we get precisely the \psd\ \ops\ of the form (\ref{psddef}).

Given a skew-symmetric $n\times n$ matrix $J$ and $F\in\bcn$, let us denote by $L_F$ the \psd\ \op\
$a(x,D)\in{\mathcal B}^*(E)$ with symbol $a(x,\xi)=F(x-J\xi)$. Let us further denote by $R_F$ the \psd\ \op\
with symbol $b(x,\xi)=F(x+J\xi)$ defined similarly as in (\ref{psddef}), except that $b(x,\xi)$
multiplies $\hat u(\xi)$ on the right. At the end of Chapter~4 in
\cite{R}, Rieffel made a conjecture that may be rephrased as follows: any
$B\in\h$ which commutes with every $R_G$, $G\in\bcn$, is of the form 
$B=L_F$ for some $F\in\bcn$. 

Using Cordes' characterization of the Heisenberg-smooth operators in the scalar
case, we have shown \cite{MM1} that Rieffel's conjecture is true when
${\mathcal A}=\C$. The second author \cite[Theorem 3.5]{M} proved further that Rieffel's
conjecture is true for any separable $C^*$-algebra ${\mathcal A}$ for which the operator 
$O_{\mathcal A}$ is a bijection. 

The assumption of separability of ${\mathcal A}$ is needed to justify several results
about vector-valued integration \cite[Ap\^endice]{Mt}.

\section{Adjointable operators}

Let us denote by $E_k$ the Hilbert module obtained using the procedure described in the first
paragraph of this paper with ${\mathcal A}$ replaced by $M_k({\mathcal A})$, the C$^*$-algebra 
of $k$-by-$k$ matrices with entries in ${\mathcal A}$.

Using that the norm $||((a_{ij}))_{1\leq i,j\leq k}||_\infty:=\max\{||a_{i,j}||;1\leq i,j\leq k\}$ is
equivalent to the C$^*$-norm $||\cdot||$ of $M_k({\mathcal A})$
($||\cdot||_\infty\leq||\cdot||\leq k^2||\cdot||_\infty$), one easily proves that a given function 
\[
f=((f_{ij}))_{1\leq i,j\leq k}:\rn\to M_k({\mathcal A})
\] 
belongs to ${\mathcal S}^{M_{k}({\mathcal A})}(\rn)$
if and only if each $f_{ij}$ belongs to ${\mathcal S}^{\mathcal A}(\rn)$.

\begin{pro} For each $(l,m)$, $1\leq l,m\leq k$, the maps
\[
P_{lm}:{\mathcal S}^{M_{k}({\mathcal A})}(\rn)\ni((f_{ij}))_{1\leq i,j\leq k}\longmapsto 
f_{lm}\in{\mathcal S}^{\mathcal A}(\rn)
\]and\[
I_{lm}:{\mathcal S}^{\mathcal A}(\rn)\ni f\longmapsto ((\delta_{il}\delta_{jm}f))_{1\leq i,j\leq k}
\in{\mathcal S}^{M_{k}({\mathcal A})}(\rn)
\]
$(\delta_{pq}=1$ if $p=q$ and $\delta_{pq}=0$ if $p\neq q)$ extend continuously to
\[
P_{lm}:E_k\longrightarrow E\ \ \mbox{and}\ \ I_{lm}:E\longrightarrow E_k.
\]
Moreover, $||P_{lm}||=1$ and $I_{lm}$ is an isometry.
\end{pro}

\pf 
For each $f\in{\mathcal S}^{\mathcal A}(\rn)$ and each $(i,j)$, we have:
\[
\left(\int I_{lm}(f)(x)^*I_{lm}(f)(x)\,dx\right)_{ij}=\delta_{im}\delta_{jm}\int f(x)^*f(x)\,dx;
\]
and, then,
\[
||I_{lm}(f)||_2^2=||((\delta_{im}\delta_{jm}\um))_{1\leq i,j\leq k}||\cdot||f||_2^2=
||f||_2^2.
\]
This shows that $I_{lm}$ is as an isometry

Given $f=((f_{ij}))_{1\leq i,j\leq k}\in{\mathcal S}^{M_{k}({\mathcal A})}(\rn)$, we have:
\[
||P_{ml}(f)||_2^2=\left|\left|\int f_{ml}(x)^*f_{ml}(x)\,dx\right|\right|\leq
\left|\left|\sum_{j=1}^{k}\int f_{jl}(x)^*f_{jl}(x)\,dx\right|\right|
\]
(we have used that $||a||\leq||a+b||$ if $a$ and $b$ are two positive elements of any C$^*$-algebra, and 
that the integral of a positive valued function is also positive). The right-hand side of the previous inequality
equals
\[
\left|\left|\left(\int f(x)^*f(x)\,dx\right)_{ll}\right|\right|\leq
\left|\left|\int f(x)^*f(x)\,dx\right|\right|_\infty\leq
\left|\left|\int f(x)^*f(x)\,dx\right|\right|.
\]
This shows that $||P_{ml}||\leq 1$. The equality holds because, for any $g\in\sC$, 
$P_{ml}(I_{ml}(g))=g$ and $||g||_2=||I_{ml}(g)||_2$.
\cqd

\begin{pro}\label{matriz}
The map
\begin{equation}\label{iso}
E_k\ni f\longmapsto ((P_{ij}(f)))_{1\leq i,j\leq k}\in\bigoplus_{1\leq i,j\leq k}E
\end{equation}
is a Banach space isomorphism. The right action of $M_k({\mathcal A})$ on $E_k$ is 
then given by matrix multiplication, while the $M_k({\mathcal A})$-valued inner-product 
on $E_k$ is given by:
\begin{equation}
\label{mvip}
 \langle\,((f_{ij}))_{1\leq i,j\leq k},((g_{ij}))_{1\leq i,j\leq k}\,\rangle=
\left(\left(\sum_{l=1}^{k}\langle f_{li},g_{lj}\rangle\right)\right)_{1\leq i,j\leq k}.
\end{equation}
\end{pro}

\pf Using that $P_{lm}I_{lm}$ equals the identity on $E$ for every $(l,m)$, 
that $P_{lm}I_{pq}=0$ if $l\neq p$ or $m\neq q$ and that $\sum_{lm}I_{lm}P_{lm}$ 
equals the identity on $E_k$, it follows that 
\[
\bigoplus_{i,j=1}^{k}E\ni((f_{ij}))_{1\leq i,j\leq k}
\longmapsto
\sum_{l,m=1}^{k}I_{lm}(f_{lm})\in E_k
\]
is an inverse for the map defined in (\ref{iso}). The statements about the action of $M_k({\mathcal A})$
and about the inner-product follow by density, since they hold on ${\mathcal S}^{M_{k}({\mathcal A})}(\rn)$.\cqd

Let ${\mathcal L}(E_k)$ denote the algebra of bounded operators on the Banach space
$E_k$. In other to describe which elements of ${\mathcal L}(E_k)$ 
belong to ${\mathcal B}^*(E_k)$ (i.e., which of them are adjointable $M_k({\mathcal A})$-module 
homomorphisms), it is convenient to define an isomorphism
\begin{equation}\label{iso2}
\bigoplus_{i,j=1}^{k}E\simeq
\bigoplus_{p=1}^{k^{2}}E,
\end{equation} 
using the bijection $\phi:\{1,\cdots,k\}\times\{1,\cdots,k\}\to \{1,2,\cdots,k^2\}$ defined by 
listing the pairs $(l,m)$ column after column,
\[
\phi(1,1)=1,\cdots,\phi(k,1)=k,\phi(1,2)=k+1,\cdots,\phi(k,2)=2k,
\]\[\cdots,
\phi(1,k)=k^2-(k-1),\cdots,
\phi(k,k)=k^2.
\]
The composition of the two isomorphisms 
defined in (\ref{iso}) and (\ref{iso2}) induces the isomorphism 
\begin{equation}
\label{iso3}
E_k\simeq\bigoplus_{p=1}^{k^{2}}E,
\end{equation}
which, by its turn, induces
\begin{equation}
\label{iso4}
{\mathcal L}(E_k)\ni T\mapsto ((P_pTI_q))_{1\leq p,q\leq k^2}\in M_{k^{2}}({\mathcal L}(E_k)).
\end{equation}
Here, abusing notation, we have written $P_p$ and $I_q$ where we really meant 
$P_{\phi^{-1}(p)}$ and $I_{\phi^{-1}(q)}$.


The following theorem is purely algebraic and could be stated for general rings and modules.

\begin{thm}\label{blocos}
Using the isomorphism $(\ref{iso4})$ as an identification, a given 
\[
T=((T_{pq}))_{1\leq p,q\leq k^2}\in{\mathcal L}(E_k)
\]
is a $($right$)$ $M_k({\mathcal A})$-module homomorphism if and only if
\begin{equation}\label{blocks}
T = \left[ \begin{array}{cccc} \tilde T & 0 & \cdots & 0 \\
0 & \tilde T & \cdots & 0 \\ \cdot & \cdot && \cdot \\ \cdot & \cdot && \cdot 
\\ \cdot & \cdot && \cdot \\ 0 & 0 & \cdots & \tilde T\end{array}\right],
\end{equation}
where $\tilde T$ is a $k$-by-$k$ matrix of bounded $($right$)$ ${\mathcal A}$-module homomorphisms and $0$ 
denotes the $k$-by-$k$ zero block.
\end{thm}
\pf Given $T=((T_{pq}))_{1\leq p,q\leq k^2}\in{\mathcal L}(E_k)$, each $T_{pq}=P_pTI_q$ is obviously bounded. If $T$ is 
an $M_k({\mathcal A})$-module homomorphism, then $T_{pq}$ is an ${\mathcal A}$-module homomorphism since, for
every $a\in{\mathcal A}$ and $f\in E$, we have 
\[
I_q(fa)=I_q(f)\left[\begin{array}{cccc}
a&0&\cdots&0\\
0&a&\cdots&0\\
\cdot&\cdot&&\cdot\\
\cdot&\cdot&&\cdot\\
\cdot&\cdot&&\cdot\\
0&0&\cdots&a\end{array}\right].
\]

Given an integer $l$, $1\leq l\leq k^2$, let $l_1$ and  $l_2$ be the integers defined by 
$0\leq l_1\leq k-1$, $1\leq l_2\leq k$ and $l=kl_1+l_2$. The product of two matrices can then be expressed by
\[
\left[
\begin{array}{cccc}
a_1&a_{1+k}&\cdots&a_{1+(k-1)k}\\
a_2&a_{2+k}&\cdots&a_{2+(k-1)k}\\
\cdot&\cdot&&\cdot\\
\cdot&\cdot&&\cdot\\
\cdot&\cdot&&\cdot\\
a_k&a_{k+k}&\cdots&a_{k+(k-1)k}\\
\end{array}
\right]
\cdot
\left[
\begin{array}{cccc}
b_1&b_{1+k}&\cdots&b_{1+(k-1)k}\\
b_2&b_{2+k}&\cdots&b_{2+(k-1)k}\\
\cdot&\cdot&&\cdot\\
\cdot&\cdot&&\cdot\\
\cdot&\cdot&&\cdot\\
b_k&b_{k+k}&\cdots&b_{k+(k-1)k}\\
\end{array}
\right]=
\]\[
=\left[
\begin{array}{cccc}
c_1&c_{1+k}&\cdots&c_{1+(k-1)k}\\
c_2&c_{2+k}&\cdots&c_{2+(k-1)k}\\
\cdot&\cdot&&\cdot\\
\cdot&\cdot&&\cdot\\
\cdot&\cdot&&\cdot\\
c_k&c_{k+k}&\cdots&c_{k+(k-1)k}\\
\end{array}
\right],
\]
with 
\[
c_l=\sum_{j=1}^{k}a_{l_2+k(j-1)}b_{j+kl_1}.
\]
This formula holds if the two matrices that we multiply belong to $M_k({\mathcal A})$ or if 
the left one is an element of $E_k$ regarded as a matrix by (\ref{iso}). With this notation, 
if a given $T=((T_{pq}))_{1\leq p,q\leq k^2}\in{\mathcal L}(E_k)$ is an $M_k({\mathcal A})$-module 
homomorphism, then, for every $(a_l)_{1\leq l\leq k^2}\in E_k$ (we now refer to 
the isomorphism (\ref{iso3})), and for every 
\[
\left[
\begin{array}{cccc}
b_1&b_{1+k}&\cdots&b_{1+(k-1)k}\\
b_2&b_{2+k}&\cdots&b_{2+(k-1)k}\\
\cdot&\cdot&&\cdot\\
\cdot&\cdot&&\cdot\\
\cdot&\cdot&&\cdot\\
b_k&b_{k+k}&\cdots&b_{k+(k-1)k}\\
\end{array}
\right]\in M_k({\mathcal A}),
\]
we have, for every $1\leq p\leq k^2$,
\begin{equation}
\label{amg}
\sum_{l=1}^{k^2}\sum_{j=1}^{k}T_{pl}a_{l_2+k(j-1)}b_{j+kl_1} =
\sum_{j=1}^{k}\sum_{l=1}^{k^2}T_{p_2+k(j-1),l}a_{l}b_{kp_1+j}.
\end{equation}
We now apply this equality, for each $a\in E$ and each $(q,r)$, $1\leq q,r\leq k^2$, to $a_l=\delta_{ql}a$
and $b_l=\delta_{rl}\um$. The only nonvanishing term in the left side sum will satisfy 
$l_2+k(j-1)=q$ (hence $j-1=q_1$ and $l_2=q_2$) and $kl_1+j=r$ (hence $l_1=r_1$ and $j=r_2$);
hence $l=kr_1+q_2$ and $j=q_1+1=r_2$. The only nonvanishing term in the right side sum will satisfy 
$l=q$ and $kp_1+j=r$ (hence $p_1=r_1$ and $j=r_2$). Equation (\ref{amg}) then becomes
\[
T_{p,kr_1+q_2}\delta_{q_1+1,r_2}=T_{p_2+k(r_2-1),q}\delta_{p_1,r_1}.
\]
We now can see that, if $q_1+1\neq r_2$ and $p_1=r_1$, then $T_{p_2+k(r_2-1),q}=0$. This proves
that, for each $(p,q)$, $T_{pq}=0$ unless $p=p_2+kq_1$. In other terms, if $p_1\neq q_1$, then 
$T_{p,q}=0$. Therefore, all blocks outside the diagonal in (\ref{blocks}) indeed vanish.

Finally, letting $q_1+1=r_2$ and $p_1=r_1$, we get $T_{p,kp_1+q_2}=T_{p_2+kq_1,q}$, or
$T_{kp_1+p_2,kp_1+q_2}=T_{kq_1+p_2,kq_1+q_2}$, proving that all blocks along the 
diagonal in (\ref{blocks}) are indeed equal. 

We have proven that any module homomorphism on $E_k$ is of the form (\ref{blocks}). To 
see that the converse is also true, we only need to remark that, under the description of $E_k$ 
given by Proposition~\ref{matriz}, 
the action of a $T$ as in (\ref{blocks}) on $E_k$ is given by left multiplication by $\tilde T$.\cqd

\begin{thm}
\label{bstar}
A given
$
T=((T_{pq}))_{1\leq p,q\leq k^2}\in{\mathcal L}(E_k)
$
belongs to  ${\mathcal B}^*(E_k)$ if and only if it is of the form $(\ref{blocks})$ with 
$\tilde T\in M_k({\mathcal B}^*(E))$.
\end{thm}

\pf Given $T$ and $S$ in ${\mathcal L}(E_k)$ of the form (\ref{blocks}), with corresponding 
$\tilde T=((\tilde T_{ij}))_{1\leq i,j\leq k}$ and 
$\tilde S=((\tilde S_{ij}))_{1\leq i,j\leq k}$, and given
$f=((f_{ij}))_{1\leq i,j\leq k}$ and  $g=((g_{ij}))_{1\leq i,j\leq k}$ in 
$E_k\simeq\bigoplus_{1\leq i,j\leq k}E$, by (\ref{mvip}) we have:
\begin{equation}\label{adjoint}
\langle Tf,g\rangle_{ij}=\sum_{l,m=1}^{k}\langle \tilde T_{lm}f_{mi},g_{lj}\rangle\ 
\mbox{and}\ \sum_{l,m=1}^{k}\langle f_{mi},\tilde S_{lm}g_{lj}\rangle=\langle f,Sg\rangle_{ij}.
\end{equation}
From this, it follows that, if each $\tilde T_{ij}$ is adjointable and if $\tilde S_{ij}=T_{ji}^*$
for all $(i,j)$, then $S$ is the adjoint of $T$. 

Conversely, suppose that $T$ is adjointable and that its adoint is $S$. The equality of the two sums
in (\ref{adjoint}) for particular choices of $f$ and $g$ will imply that each $\tilde T_{ij}$ 
is adjointable. \cqd

\section{Heisenberg-smooth adjointable operators}
The mapping (\ref{azz}) may be defined for any $B$ in ${\mathcal L}(E)$ or in ${\mathcal L}(E_k)$. 
It thus makes sense to talk about Heisenberg-smooth operators in ${\mathcal L}(E)$ or in 
${\mathcal L}(E_k)$. Given $B=((B_{pq}))_{1\leq p,q\leq k^2}\in{\mathcal L}(E_k)$, we have
\[
B_{z,\zeta}=(((B_{pq})_{z,\zeta}))_{1\leq p,q\leq k^2}=((P_pB_{z,\zeta}I_q))_{1\leq p,q\leq k^2}
\]
(it is enough to check these equalities on the dense subset of Schwartz functions).
We then get:

\begin{pro}\label{cada}
A given  $B=((B_{pq}))_{1\leq p,q\leq k^2}\in{\mathcal L}(E_k)$ is Heisenberg-smooth if and only if
each $B_{pq}\in{\mathcal L}(E)$ is Heisenberg-smooth. 
\end{pro}

This leads to our main result:
\begin{thm}
If the unital separable C$^*$-algebra ${\mathcal A}$ is such that the map $O_{{\mathcal A}}$ defined
in $(\ref{o})$ is a bijection, then the map $O_{M_k({\mathcal A})}$ is also a bijection.
\end{thm}
\pf
Given any Heisenberg-smooth operator $B\in{\mathcal B}^*(E_k)$, we have to show that it is of the form 
$B=b(x,D)$ for some $b\in\bCk$. Let $((B_{pq}))_{1\leq p,q\leq k^2}$ be the matrix that corresponds to 
$B$ by the isomorphism (\ref{iso4}). By Proposition~\ref{cada} and
by the assumption that $O_{{\mathcal A}}$ is surjective, each $B_{pq}$ is a pseudodifferential operator
of the type defined in (\ref{psddef}). By Theorem~\ref{blocos}, $B$ is of the form (\ref{blocks}). That is, 
there exist $b_{ij}\in\bC$, $1\leq i,j\leq k$, such that, with
\[ 
\tilde T = \left[ \begin{array}{cccc} b_{11}(x,D) & b_{12}(x,D) & \cdots & b_{1k}(x,D) \\
b_{21}(x,D) &b_{22}(x,D) & \cdots & b_{2k}(x,D) \\ \cdot &\cdot&& \cdot \\ 
\cdot &\cdot&& \cdot \\ \cdot &\cdot&& \cdot 
\\b_{k1}(x,D) &b_{k2}(x,D) & \cdots & b_{kk}(x,D)  \end{array}\right],
\]
we have
\[ 
B = \left[ \begin{array}{cccc} \tilde T & 0 & \cdots & 0 \\
0 & \tilde T & \cdots & 0 \\ \cdot &\cdot&& \cdot \\ \cdot &\cdot&& \cdot \\ 
\cdot &\cdot&& \cdot \\ 0 & 0 & \cdots & \tilde T\end{array}\right].
\]
This implies that $B$ and $b(x,D)$ are equal, if $b\in\bCk=M_k(\bC)$ is given by 
$b=((b_{ij}))_{1\leq i,j\leq k}$. Indeed, the equality of the two \ops\ can be easily 
verified on ${\mathcal S}^{M_{k}({\mathcal A})}(\rn)$. \cqd

\section*{Acknowledgements} Severino Melo was partially supported by the
Brazilian agency CNPq (Processo 306214/2003-2), and Marcela Merklen had a posdoc
position sponsored by Fapesp (Processo 2006/07163-9). We thank 
Fernando Abadie and H\'ector Merklen for some very helpful conversations.

\vskip1.0cm

{\footnotesize 

\noindent
Instituto de  Matem\'atica e Estat\'{\i}stica, Universidade de S\~ao Paulo,\\
Rua do Mat\~ao 1010, 05508-090 S\~ao Paulo, Brazil. 

\noindent
Email: toscano@ime.usp.br, marcela@ime.usp.br}

\end{document}